\documentclass[10pt, a4paper, fleqn]{scrartcl}
\usepackage[latin1]{inputenc}
\usepackage{a4wide}
\usepackage{amsmath}
\usepackage{amssymb}
\usepackage[all]{xy}
\usepackage{amsthm}
\usepackage{paralist, calc}
\usepackage{graphicx}

\usepackage{color} % Colored text

\newcommand\F{\ensuremath{\mathbb{F}}}

\newcommand\Fq{\ensuremath{\mathbb{\F}_{q^2}}}

\renewcommand\phi{\ensuremath{\varphi}}
\renewcommand\epsilon{\ensuremath{\varepsilon}}

\newcommand\mrm[1]{\mathrm{#1}}

\newtheorem{thm}{Theorem}[section]
\newtheorem{cor}[thm]{Corollary}
\newtheorem{lemma}[thm]{Lemma}
\newtheorem{prop}[thm]{Proposition}

\theoremstyle{remark}
\newtheorem{rem}[thm]{Remark}

\theoremstyle{definition}
\newtheorem{defi}[thm]{Definition}

\numberwithin{equation}{section}

\newenvironment{Proof}[1]{\noindent {\em Proof #1.}}{\hfill $\Box$ \medskip}

\numberwithin{equation}{section}

\begin{document}

\title{{\bf Isomorphisms of unitary forms of Kac-Moody groups over finite fields}}
\author{Ralf Gramlich\thanks{The first author gratefully acknowledges a Heisenberg fellowship by the Deutsche Forschungsgemeinschaft.} \and Andreas Mars}
\maketitle

\begin{abstract}
We use the methods developed in \cite{main-caprace}, \cite{caprace-muehlherr-isomorphisms}, \cite{caprace-muehlherr-isomorphisms-bounded} to solve the isomorphism problem of unitary forms of infinite split Kac-Moody groups over finite fields of square order.  
\end{abstract}

\section{Introduction}

The isomorphism problem for Kac-Moody groups has been studied and solved in \cite{main-caprace}, \cite{caprace-muehlherr-isomorphisms}, \cite{caprace-muehlherr-isomorphisms-bounded}. In addition, \cite{main-caprace} contains a description of the isomorphisms of unitary forms of complex Kac-Moody groups with respect to the compact involution; see also \cite{kac-peterson-defining-relations}.

An adaption of the methods developed in \cite{main-caprace}, \cite{caprace-muehlherr-isomorphisms}, \cite{caprace-muehlherr-isomorphisms-bounded} allows us to prove the following theorem. (We refer to Section \ref{basics} for definitions.)

\medskip \noindent {\bf Main Result.} 
{\em 
Let $q$ and $r$ be arbitrary prime powers, let $G$ and $G'$ be infinite split Kac-Moody groups over $\mathbb{F}_{q^2}$ and $\mathbb{F}_{r^2}$, respectively, and let $K$ and $K'$ be their respective unitary forms.

If there exists an isomorphism $\phi: K \rightarrow K'$, then the following hold.
\begin{enumerate}
\item $q = r$.
\item There exist a bijection $\pi: S \rightarrow S'$, an inner automorphism $\nu$ of $K'$ and for each $i \in S$ a
diagonal-by-field automorphism $\gamma_i$ of $\mrm{SU}_2(\Fq)$ such that the diagram
  \[
     \xymatrix{
        \mrm{SU}_2(\Fq)\ar[dd]_{\phi_i} \ar[rr]^{\gamma_i} & & \mrm{SU}_2(\Fq)\ar[dd]^{\phi_{\pi(i)}'} \\
        & & \\
        K\ar[rr]^{\nu \circ \phi} & & K'
     }
  \]
commutes for every $i \in S$.
\end{enumerate}
}

As a byproduct of our strategy, one may conclude that the underlying Kac-Moody root data are isomorphic (cf.\ \cite[Section 2.3.2]{caprace-muehlherr-isomorphisms-bounded}). This implies the following, because $N_G(K) = K$. (In fact, for sufficiently large $q$ one can deduce from \cite[Corollary 7.7(ii)]{Caprace/Monod} and \cite{Gramlich/Muehlherr} that $\mathrm{Comm}_G(K) = K$.)

\medskip \noindent {\bf Corollary (Strong rigidity).}
{\em 
For any isomorphism $\phi : K \rightarrow K'$ there exists a unique isomorphism $\psi : G \to G'$ satisfying $\psi|_K^{K'} \equiv \phi$. 
}

\medskip
Throughout the paper we assume that the reader is familiar with Kac-Moody groups and their buildings as described, for instance, in \cite{Abramenko/Brown}, \cite[Chapter 1]{main-caprace}. Section \ref{basics} may serve as a reminder about these concepts. In Section \ref{finite-subgroups} we collect information about the structure of maximal finite subgroups, and in Section \ref{isomorphisms} we prove the Main Result.

\medskip
This note can be considered as being part of an ongoing project of understanding unitary forms of Kac-Moody groups over finite fields. We are particularly interested in these unitary forms as over sufficiently large finite fields they are lattices in the completions of the ambient Kac-Moody groups with respect to the topology of compact convergence (cf.\ \cite{Gramlich/Muehlherr}).

\subsubsection*{Acknowledgements}
The authors express their gratitude to Max Horn and Stefan Witzel for several helpful discussions on the topic of this article. Moreover, the authors thank an anonymous referee for an excellent report with a lot of detailed comments and suggestions. Finally, we gratefully acknowledge that we owe the strategy of the proof presented in this final version to Pierre-Emmanuel Caprace (personal communication) and to the referee.

%%%%%%%%%%%%%%%%%%%%%%%
%%%%%%%%%%%%%%%%%%%%%%%
%%% SECTION: Basics %%%
%%%%%%%%%%%%%%%%%%%%%%%
%%%%%%%%%%%%%%%%%%%%%%%

\section{Basics} \label{basics}
Let $G$ be a split Kac-Moody group over a finite field (cf.\ \cite[Section 8.11]{Abramenko/Brown}, \cite[Chapter 3]{main-caprace}, \cite[Chapter 8]{remy-kac-moody}, \cite{tits-algebra}) and let $(G,B_+,B_-,N,S)$ be the associated saturated twin Tits system with Weyl group $W$. For $\epsilon \in \{+,-\}$ we have the \textit{Bruhat}, resp.\ \textit{Birkhoff decompositions}
  \[
    G = \bigsqcup_{w\in W} B_\epsilon wB_\epsilon \quad \quad \quad \quad G = \bigsqcup_{w\in W} B_\epsilon w B_{-\epsilon},
  \]
see \cite[Section 6.2 and Proposition 6.81]{Abramenko/Brown}, \cite[Chapter 1]{main-caprace}, \cite[Chapter 1]{remy-kac-moody}.

Conjugates of the {\em fundamental Borel subgroups} $B_+$ and $B_-$ are called \textit{Borel subgroups} of $G$. The intersection $T := B_+ \cap B_-$ is called the {\em fundamental maximal split torus} of $G$; each of its conjugates is called a {\em maximal split torus}. A \textit{fundamental parabolic subgroup} $P_\epsilon$ of $G$ is a subgroup containing a fundamental Borel group $B_\epsilon$. Any conjugate of a fundamental parabolic subgroup is simply called \textit{parabolic subgroup}. To a fundamental parabolic subgroup $P_\epsilon$ of $G$ there exists $J \subseteq S$ such that
\[
	P = \bigsqcup_{w\in W_J} B_\epsilon wB_\epsilon
\]
for the {\em special subgroup} $W_J := \langle J \rangle$ of $W$. The set $J$ is called the \textit{type} of $P$. Moreover, $P$ and $J$ are called \textit{spherical}, if $W_J$ is finite.

\medskip
A split Kac-Moody group $G$ can be defined functorially, cf.\ \cite[Chapter 8]{remy-kac-moody}, \cite{tits-algebra}. In particular, for each $s \in S$ there exists a homomorphism $\phi_s : \mathrm{SL}_2(\Fq) \to G$ with central kernel such that $$G = \langle \phi_s(\mathrm{SL}_2(\Fq)) \mid s \in S \rangle.$$ 

\begin{defi}[Unitary form]
Let $G$ be a split Kac-Moody group over the field $\mathbb{F}_{q^2}$, let $\omega$ be the Chevalley involution of $G$, cf.\ \cite[Chapter 8]{main-caprace}, \cite[Section 2]{kac-peterson-defining-relations}, and let $\theta$ be the composition of $\omega$ and the field involution $x \mapsto x^q$ of $\mathbb{F}_{q^2}$, called \textit{twisted Chevalley involution}. The fixed point group $K := \{ g \in G \mid \theta(g) = g \}$ is called the \textit{unitary form} of $G$ with respect to $\theta$. 
\end{defi}

This unitary form is the finite field analogue of the real group studied in \cite[Section 5]{kac-peterson-defining-relations}, where the involved field involution is complex conjugation.

\medskip
For each $s \in S$, the intersection $K_s := \phi_s(\mathrm{SL}_2(\mathbb{F}_{q^2})) \cap K$ is isomorphic to $\mathrm{SU}_2(\mathbb{F}_{q^2})$ and is called a \textit{rank one subgroup} of $K$. The intersection $T_K := T \cap K$ is called the {\em fundamental torus} of $K$. The involution $\theta$ induces an involution of $\phi_s(\mathrm{SL}_2(\mathbb{F}_{q^2}))$ which pulls back to the product of the contragredient automorphism of $\mathrm{SL}_2(\mathbb{F}_{q^2})$ and the field involution, whose set of fixed elements forms a subgroup isomorphic to $\mathrm{SU}_2(\mathbb{F}_{q^2})$.

\begin{defi}[Twin building] \label{twin}
Let $G$ be a split Kac-Moody group and let $(G,B_+,B_-,N,S)$ be the associated saturated twin Tits system with Weyl group $W$.
Let $\Delta_\epsilon := G/B_\epsilon$, let
\begin{eqnarray*}
    \delta_\epsilon: \Delta_\epsilon \times \Delta_\epsilon & \to & W \\
  \delta_\epsilon(gB_\epsilon,hB_\epsilon) & := & w \quad \mbox{ if and only if $B_\epsilon g^{-1}hB_\epsilon = B_\epsilon w B_\epsilon$,}
\end{eqnarray*}
and let
\begin{eqnarray*}
\delta_*: (\Delta_+ \times \Delta_-) \cup (\Delta_- \times \Delta_+) & \to & W \\
    \delta_*(gB_\epsilon,hB_{-\epsilon}) & := & w \quad \mbox{ if and only if $B_\epsilon g^{-1}h B_{-\epsilon} = B_\epsilon w B_{-\epsilon}$.}
\end{eqnarray*}
The triple $((\Delta_+,\delta_+),(\Delta_-,\delta_-),\delta_*)$ is the {\em twin building} associated to $G$, cf.\ \cite[Chapter~8]{Abramenko/Brown}.
\end{defi}

Note that the maps $\delta_\epsilon$ and $\delta_*$ are well-defined, because of the Bruhat and Birkhoff decompositions. 

\medskip
Each pair $(\Delta_\epsilon,\delta_\epsilon)$ is a building. The {\em Davis realisation} of a building, as described in \cite{Davis:1997}, is a $\mathrm{CAT}(0)$ space. In particular, the Bruhat-Tits fixed point theorem can be applied to the Davis realisation. For details on $\mathrm{CAT}(0)$ spaces, we refer the reader to the book \cite{bridson-haefliger}.

\begin{defi}[Flipflop system]
Let $G$ be a split Kac-Moody group over $\mathbb{F}_{q^2}$ and let $\theta$ be the twisted Chevalley involution of $G$. The involution $\theta$ induces an involutory automorphism of the twin building associated to $G$, which we also denote by $\theta$.
A \textit{Phan chamber} is a chamber $c \in \Delta$ satisfying $\delta_*(c, \theta(c)) = 1_W$. The collection of all Phan chambers contained in $(\Delta_+,\delta_+)$ is called the \textit{flipflop system} and is denoted by $\Delta_\theta$.
\end{defi}

%%%%%%%%%%%%%%%%%%%%%%%%%%%%%%%%%%%%%%%%%%%%%%%%%%%%%%%%%%%%%%
%%%%%%%%%%%%%%%%%%%%%%%%%%%%%%%%%%%%%%%%%%%%%%%%%%%%%%%%%%%%%%
%%% SECTION 3 about finite subgroups                       %%%
%%%%%%%%%%%%%%%%%%%%%%%%%%%%%%%%%%%%%%%%%%%%%%%%%%%%%%%%%%%%%%
%%%%%%%%%%%%%%%%%%%%%%%%%%%%%%%%%%%%%%%%%%%%%%%%%%%%%%%%%%%%%%

\section{On maximal finite subgroups} \label{finite-subgroups}
Our strategy is to show that any isomorphism between two unitary forms induces a bijection between the respective sets of rank one subgroups. This implies that an isomorphism is uniquely determined by its behaviour on the rank one subgroups.

We first investigate the structure of maximal finite subgroups.

%%%%%%%%%%%%%%%%%%%%%%%%%%%%%%%%%%%%%%%%%%%%%%%%%%%%%%%%%%%%%%%%%%%%%%%%%%%%%%%%%%%%%%%%%%%%%%%%%%%%%
%%% Proposition: \theta-fixed point in unipotent radical of spherical parabolics of opposite sign %%%
%%%%%%%%%%%%%%%%%%%%%%%%%%%%%%%%%%%%%%%%%%%%%%%%%%%%%%%%%%%%%%%%%%%%%%%%%%%%%%%%%%%%%%%%%%%%%%%%%%%%%
\begin{prop} \label{intersections}
Let $q$ be a prime power.  Let $G$ be a split Kac-Moody group over $\mathbb{F}_{q^2}$ and let $K = \rm{Fix}_G(\theta)$ be the unitary form of $G$ with respect to the twisted Chevalley involution $\theta$.

Let $P_1, P_2 \leq G$ be spherical parabolic subgroups of opposite sign of $G$ and set
\[
	P := P_1 \cap P_2 = L(P) \ltimes U(P),
\]
where $U(P)$ is the unipotent radical of $P_1 \cap P_2$. If $U(P)$ is non-trivial and $\theta$ normalises $P_1 \cap P_2$, then $K \cap U(P)$ is non-trivial.
\end{prop}

%%% Remark: About Levi decomposition etc. %%%
\begin{rem} \label{rem:prop3.1}
\begin{enumerate}
\item Let $p = \mathrm{char}(\Fq)$. The above decomposition is the Levi decomposition
\[
	P = L(P) \ltimes U(P),
\]
given in \cite[Proposition 3.6]{caprace-muehlherr-isomorphisms-bounded}, i.e.\ $L(P)$ is a semisimple Levi subgroup of $P$ and $U := U(P)$ is a $p$-group. 

\item Note that since $P_1$ and $P_2$ are spherical parabolic subgroups of opposite sign, by definition $P_1 \cap P_2$ is bounded. However, for Kac-Moody groups over finite fields the terms finite and bounded are equivalent, see \cite[Lemma 6.2]{main-caprace}. Therefore $P$ is a finite group.
%
%\item Moreover $U$ is the unipotent radical of $P$ which is generated by root groups as described in \cite[Section 3.2.1]{caprace-muehlherr-isomorphisms-bounded}; however, we will not need the precise description here. For $q \in \{2,3\}$, we have to add the unipotent radicals of the fundamental Borel groups of the fundamental rank two groups as generators for $U$ in order to obtain a finite set of generators.
\end{enumerate}
\end{rem}

\begin{proof}[Proof of Proposition \ref{intersections}.]
Let $R_{P_+}$ and $R_{P_-}$ be the maximal spherical residues of $\Delta_+$, respectively $\Delta_-$ (see Definition \ref{twin}), that are stabilised by $P_+$ and $P_-$ respectively. Since $P_- = \theta(P_+)$, there exists a twin apartment $\Sigma$ satisfying $\theta(\Sigma \cap R_{P_+}) = \Sigma \cap R_{P_-}$. We can assume that $R_{P_+}$ and $R_{P_-}$ are not opposite, as otherwise $U(P)$ would be trivial, and there were nothing to prove. Therefore, using the construction of $\theta$-stable apartments in \cite[Chapter 2]{Horn:2008}, there exist two distinct $\theta$-stable twin apartments $\Sigma'$ and $\Sigma''$ both containing $\Sigma \cap R_{P_+}$ and $\Sigma \cap R_{P_-}$, because $\theta$ involves the field involution of $\mathbb{F}_{q^2}$. As $\Sigma$' and $\Sigma''$ are conjugate  by an element $g \in G$ that lies both in $K$ (see \cite[Chapter 2]{Horn:2008}) and in $U(P)$ (see \cite[Proposition 3.1]{caprace-muehlherr-isomorphisms-bounded}), the intersection $K \cap U(P)$ is necessarily non-trivial.
\end{proof}

%%%%%%%%%%%%%%%%%%%%%%%%%%%%%%%%%%%%%%%%%%%%%%%%%%%%%%%%%
%%% Lemma: unitary form transitive on flipflop-system %%%
%%%%%%%%%%%%%%%%%%%%%%%%%%%%%%%%%%%%%%%%%%%%%%%%%%%%%%%%%
% [AM] Need three-spherical here according to Theorem 6.6 of "simple connectedness of certain subsets of buildings", right?
\begin{lemma} \label{lem:unitary-form-transitive-on-ffs}
Let $K$ be the unitary form of a split Kac-Moody group $G$ with respect to $\theta$. Then $K$ acts transitively on the flipflop system $\Delta_\theta$ of $G$.
\end{lemma}
\begin{proof}
This is a well known fact. It can be proved by a local analysis as conducted in \cite{Devillers/Muehlherr}, \cite{Horn:2008} or by a standard argument based on Lang's Theorem.
\end{proof}

Recall that the $p$-core $O_p(X)$ of a finite group $X$ is the largest normal $p$-subgroup of $X$.

%%%%%%%%%%%%%%%%%%%%%%%%%%%%%%%%%%%%%%%%%%%%%%%%%%%%%%%%%%%%%%%%%
%%% Proposition: Characterisation of maximal finite subgroups %%%
%%%%%%%%%%%%%%%%%%%%%%%%%%%%%%%%%%%%%%%%%%%%%%%%%%%%%%%%%%%%%%%%%
\begin{prop} \label{max-finite}
Let $q$ be a power of the prime $p$, let $G$ be an infinite split Kac-Moody group over $\mathbb{F}_{q^2}$, let $K$ be its unitary form, and let $P_+, P_- \leq G$ be opposite maximal spherical parabolic subgroups with the property that $\theta(P_+) = P_-$. 

Then $P_+ \cap P_- \cap K = \mathrm{Fix}_{P_+ \cap P_-}(\theta)$ is a maximal finite subgroup of $K$ with trivial $p$-core.
Conversely, any maximal finite subgroup with trivial $p$-core is obtained in this fashion.
\end{prop}
\begin{proof}
Let $R_{P_+}$ and $R_{P_-}$ be the respective residues of $\Delta_+$ and $\Delta_-$ associated to $P_+$ and $P_-$. As $P_+$ and $P_-$ are opposite, so are  $R_{P_+}$ and $R_{P_-}$, and the intersection $P := P_+ \cap P_-$ is finite and semisimple. Hence $$(R_{P_+}, R_{P_-}, {\delta_*}_{|R_{P_\pm} \times R_{P_{\mp}}})$$ is a spherical twin building. By \cite{Tits:1992}, this twin building can be canonically identified with the spherical building $\Delta(P)$ of the finite semisimple group $P = P_+ \cap P_-$. The product $\theta \circ \mathrm{proj}_{R_{P_+}}$ yields an involution on $\Delta(P)$, which we also denote by $\theta$. Recall that an element of $P$ is semisimple if and only if it is $p$-regular, and that it is unipotent if and only if it is $p$-singular. Hence $P$, and therefore $F := P \cap K = \mathrm{Fix}_P(\theta)$, have trivial $p$-cores.

Again we denote by $\Delta_\theta$ the flipflop system of $\Delta(P)$. Let $c \in \Delta_\theta$. Then $\delta(c, \theta(c)) = w_0 \in W_P$, where $w_0$ is the longest word in the Coxeter group $W_P$ associated to $P$. Hence $\Sigma := \mathrm{conv}(c, \theta(c))$ is a $\theta$-stable apartment containing two opposite Phan chambers, which implies that $\Sigma$ consists of Phan chambers only. By Lemma \ref{lem:unitary-form-transitive-on-ffs} the group $F$ acts transitively on the set of Phan chambers. Hence the orbit of $c$ under the action of $F$ on $\Delta(P)$ meets every double coset $BwB$, where $w \in W_P$ and $B := P_c$. Hence any parabolic subgroup of $P$ containing $F$ necessarily contains all double cosets of $P$ modulo $B$.
We conclude that $F$ cannot stabilise a proper residue of the building $\Delta(P)$.

The Davis realisation (see \cite{Davis:1997}) of each half $\Delta_\epsilon$ of the twin building $\Delta(G)$ admits one obvious fixed point of $F$, namely $R_{P_\epsilon}$. We claim that these are unique. To show this, suppose there is some other spherical residue $R_{Q_\epsilon}$ in $\Delta_\epsilon(G)$ which is stabilised by $F$. By maximality of $R_{P_\epsilon}$, the residue $R_{Q_\epsilon}$ cannot contain $R_{P_\epsilon}$ properly. Moreover, the residues $R_{P_\epsilon}$ and $R_{Q_\epsilon}$ must be disjoint, else their intersection would yield a proper residue of $R_{P_\epsilon}$ stabilised by $F$, a contradiction to what we established above. Now consider the projection of $R_{Q_\epsilon}$ onto $R_{P_\epsilon}$. If $F$ stabilises $R_{Q_\epsilon}$, then it also stabilises the image of the projection. Hence the projection must be surjective. In view of \cite[Corollary 2.8]{caprace-muehlherr-isomorphisms-bounded}, the residues $R_{P_\epsilon}$ and $R_{Q_\epsilon}$ are therefore opposite, which is absurd as $G$ is infinite, whence $\Delta_\epsilon(G)$ non-spherical. Hence $R_{P_\epsilon}$ is the unique fixed point of $F$ in the Davis realisation of $\Delta_\epsilon(G)$.

However, the stabiliser of $R_{P_\epsilon}$ in $G$ is $P_\epsilon$, whence the stabiliser of the residues $R_{P_+}$ and $R_{P_-}$ in $K$ is equal to $P_+ \cap P_- \cap K = P \cap K = F$. Since by the Bruhat-Tits fixed point theorem every finite subgroup fixes a point in the Davis realisation (as it is $\mathrm{CAT}(0)$), we conclude that $F$ must be maximal among finite subgroups of $K$ with trivial $p$-core.

Conversely, let $F \leq K$ be maximal finite with trivial $p$-core. Since $F$ is finite, it has a bounded orbit on $\Delta_+$, whence the Bruhat-Tits fixed point theorem again implies that there exists some positive spherical parabolic subgroup $P_+ \leq G$ with $F \leq P_+$. Without loss of generality, we may choose $P_+$ with this property and of minimal rank. Since $F$ is pointwise fixed by $\theta$, it follows from the construction that $F$ is also contained in $\theta(P_+) =: P_-$. Now $P_+$ and $P_-$ are spherical, thus $P := P_+ \cap P_-$ is a finite group. In particular, we see that $F \subseteq P \cap K$ and $P \cap K$ is finite. By maximality of $F$, we get that $F = P \cap K$.

It remains to show that the groups $P_+$ and $P_-$ constructed above are opposite and maximal spherical. The group $F$ has trivial $p$-core by assumption, hence we may apply Proposition \ref{intersections} and see that $P$ has trivial unipotent radical. Comparing this with \cite[Proposition 3.1 and Proposition 3.6]{caprace-muehlherr-isomorphisms-bounded} and the fact that $P_+$ and $P_-$ are of minimal rank, we see that $P_+$ and $P_-$ must be opposite, otherwise the unipotent radical would be non-trivial. Since $F$ is a maximal finite subgroup, the parabolic subgroups $P_+$ and $P_-$ are maximal spherical.
\end{proof}

We conclude this section by recording a first structural property of isomorphisms of unitary forms.

%%%%%%%%%%%%%%%%%%%%%%%%%%%%%%%%%%%%%%%%%%%%%%%%%%%%%%%%%%%%%%%%%%%%%%%%%%%%%%
%%% Corollary: isomorphism preserves maximal spherical Levi subgroups of K %%%
%%%%%%%%%%%%%%%%%%%%%%%%%%%%%%%%%%%%%%%%%%%%%%%%%%%%%%%%%%%%%%%%%%%%%%%%%%%%%%
\begin{cor} \label{iso-preserves-levi}
Let $\phi: K \to K'$ be an isomorphism of unitary forms. Then $\phi$ maps maximal spherical Levi subgroups of $K$ to maximal spherical Levi subgroups of $K'$.
\end{cor}
\begin{proof}
By Proposition \ref{max-finite}, a maximal spherical Levi subgroup $L$ of $K$ is associated to some reductive group $P$ such that $L = P \cap K$, where $P$ is the intersection of two opposite maximal spherical parabolic subgroups. As the image of $L$ under $\phi$ is maximal finite in $K'$, it holds (again by Proposition \ref{max-finite}) that $\phi(L) = L' \cap K' = P' \cap K'$, where $P'$ is the intersection of some opposite maximal spherical parabolic subgroups of $G'$ and $L'$ is the Levi subgroup of their intersection. Hence the image of a maximal spherical Levi subgroup is again a maximal spherical Levi subgroup.
\end{proof}

%%%%%%%%%%%%%%%%%%%%%%%%%%%%%%%%%%%%%%%%%%%%%%%%%%%%%%%%%%%%%%%%
%%%%%%%%%%%%%%%%%%%%%%%%%%%%%%%%%%%%%%%%%%%%%%%%%%%%%%%%%%%%%%%%
%%% SECTION 4 about Isomorphisms                             %%%
%%%%%%%%%%%%%%%%%%%%%%%%%%%%%%%%%%%%%%%%%%%%%%%%%%%%%%%%%%%%%%%%
%%%%%%%%%%%%%%%%%%%%%%%%%%%%%%%%%%%%%%%%%%%%%%%%%%%%%%%%%%%%%%%%

\section{Isomorphisms} \label{isomorphisms}
We continue to denote by $G$ and $G'$ split Kac-Moody groups over the fields $\Fq$ and $\F_{r^2}$, respectively, and their respective unitary forms by $K$ and $K'$.

Firstly, we record that we may recognise the characteristic of the ground field.

%%%%%%%%%%%%%%%%%%%%%%%%%%%%%%%%%%%%%%%%%%%%%%%%%%%%%%%%%%%%%%%%%%%%%%%%%%%
%%% Proposition: Characteristic encoded in orders of finite p-subgroups %%%
%%%%%%%%%%%%%%%%%%%%%%%%%%%%%%%%%%%%%%%%%%%%%%%%%%%%%%%%%%%%%%%%%%%%%%%%%%%
\begin{prop} \label{char}
Let $K$ be an infinite unitary form and let $p$ be a prime. Then the set of orders of finite $p$-subgroups of $K$ is unbounded if and only if $p = \mathrm{char}(\Fq)$.
\end{prop}
\begin{proof}
If $p \neq \mathrm{char}(\Fq)$, then the set of orders of finite $p$-subgroups of $K$ is bounded by \cite[Proposition 6.2]{caprace-muehlherr-isomorphisms-bounded}. For $p = \mathrm{char}(\Fq)$ it follows from the proof of \cite[Theorem 1]{Gramlich/Muehlherr} that the sets of orders of the stabilisers $\mathrm{Stab}_K(c)$ and of their $p$-Sylow subgroups, respectively, are unbounded.
\end{proof}

The following result is an adaption of \cite[Theorem 5.1]{caprace-muehlherr-isomorphisms-bounded} to unitary forms.

%%%%%%%%%%%%%%%%%%%%%%%%%%%%%%%%%%%%%%%%%%%%%%%%%%%%%%%%%%%
%%% Proposition: Preservation of the rank one subgroups %%%
%%%%%%%%%%%%%%%%%%%%%%%%%%%%%%%%%%%%%%%%%%%%%%%%%%%%%%%%%%%
\begin{prop} \label{preservation}
Let $K$ and $K'$ be unitary forms of infinite split Kac-Moody groups over fields $\Fq$ and $\F_{r^2}$. If there exists an isomorphism $\phi: K \to K'$, then there exists $g \in K'$ such that
\begin{enumerate}
\item $q = r$, 
\item the tori $\phi(T_K)$ and $T'_{K'}$ are conjugate under $g$, and 
\item $\{ g\phi(K_\alpha)g^{-1} \mid \alpha \in \Phi \} = \{ K'_\alpha \mid \alpha \in \Phi' \}$.
\end{enumerate}
\end{prop}
\begin{proof}
Let $p := \mathrm{char}(\mathbb{F}_{q^2}) = \mathrm{char}(\mathbb{F}_{r^2})$ (Proposition \ref{char}) and let $F$ be a maximal finite subgroup of $K$ with trivial $p$-core. Proposition \ref{max-finite} implies $F = P_+ \cap P_- \cap K$ for some opposite maximal spherical parabolic subgroups $P_+$ and $P_-$ of $G$ satisfying $\theta(P_+) = P_-$. For the same reason we can write $\phi(F)$ as $P_+' \cap P_-' \cap K'$. Define $P := P_+ \cap P_-$ and $P' := P'_+ \cap P'_-$. By the Levi decomposition (cf.\ \cite[Proposition 3.6]{caprace-muehlherr-isomorphisms-bounded}; see also Remark \ref{rem:prop3.1} of this article) the groups $P$ and $P'$ are semisimple finite groups of Lie type. Since $F$ and $\phi(F)$ are isomorphic as abstract groups and since, moreover, $F$ and $\phi(F)$ are twisted finite groups of Lie type in identical characteristics (again Proposition \ref{char}) embedded in $P$ and $P'$, respectively (Proposition \ref{max-finite}), we conclude from \cite[Table 2.2 in Section 2.2]{Gorenstein/Lyons/Solomon:1998}, \cite[Table 6 in Chapter 4]{atlas} that $q=r$ and that the buildings and the diagrams of $P$ and $P'$ coincide. 
Hence $\phi$ induces an isomorphism $\psi$ between $P$ and $P'$. Therefore, by \cite[Theorem 7.1]{caprace-muehlherr-isomorphisms-bounded}, the map $\psi$ induces an isomorphism of the twin root datum of $P$ onto the twin root datum of $P'$. In particular, $\psi$ maps rank one subgroups of $P$ to rank one subgroups of $P'$.

% [AM] Note: |T| = (q+1)^|S| only if the Kac-Moody root datum is simply connected, which is not the case according to our assumptions... However, we only need equal order of the tori here...
Let $H_K := H \cap K $ be a $K$-conjugate of the fundamental maximal split torus $T_K$ of $K$ that is contained in $P$, which is possible by Lemma \ref{lem:unitary-form-transitive-on-ffs}, because $P_+$ and $P_-$ are opposite.  The torus $H_K$ stabilises a unique twin apartment $\Sigma$ (cf.\ \cite[Definition 5.171]{Abramenko/Brown}) of the twin building of $G$, i.e., $H_K=\mathrm{Fix}_K(\Sigma)$. By the above discussion and the fact that the fundamental tori $T_K$ and $T'_{K'}$ have the same cardinality, there exists a unique twin apartment $\Sigma'_0$ of the twin building of $G'$ such that $\phi(H_K) = \mathrm{Fix}_{K'}(\Sigma'_0) \leq P'$. By Lemma \ref{lem:unitary-form-transitive-on-ffs} the twin apartment $\Sigma'_0$ is in the $K'$-orbit of the fundamental twin apartment $\Sigma'$, so that $\phi(H_K)$ is a $K'$-conjugate of the fundamental torus $T'_{K'}$ of $K'$. Hence there exists $g \in K'$ such that $\phi' := (c_g \circ \phi) (T_K) = T'_{K'}$. 
%More precisely, $g \in K'$ is some element mapping the fundamental apartment onto $\Sigma'$ via the action of $K'$.

By the same arguments together with the facts that $T_K$ and $T'_{K'}$ fix a unique apartment in the respective buildings and each rank one subgroup is contained in a maximal spherical parabolic subgroup of $K$ or $K'$, respectively, we conclude that
\[
	\{ g\phi(K_\alpha)g^{-1} \mid \alpha \in \Phi \} = \{ K'_\alpha \mid \alpha \in \Phi' \}.
\]
%and hence $\phi'$ maps the subgroups $G_\alpha \cap K \cong \mathrm{SU}_2(\Fq)$ in a one-to-one correspondence to subgroups of the form $G'_\alpha \cap K' \cong \mathrm{SU}_2(\Fq)$.
\end{proof}

We now have everything at hand to prove the main result.

%%%%%%%%%%%%%%%%%%%%%%%%%%%%%
%%% Proof of Main Theorem %%%
%%%%%%%%%%%%%%%%%%%%%%%%%%%%%
\bigskip
\begin{Proof}{of the Main Result}
Assertation (i) equals assertion (i) of Proposition \ref{preservation}.

The inner automorphism and the bijection of the index sets in (ii) are provided by assertions (ii) and (iii) of Proposition \ref{preservation}. It thus remains to analyse the isomorphism $\phi$ restricted to the rank one subgroups $K_\alpha \cong \mathrm{SU}_2(\Fq) \cong \mathrm{SL}_2(\mathbb{F}_{q})$. By \cite{steinberg-automorphisms} the outer automorphisms of $K_\alpha$ are diagonal-by-field and the result follows.
\end{Proof}

\vspace{1cm}

\noindent Authors' address:

\vspace{.8cm}

\noindent 
Ralf Gramlich, Andreas Mars \\
TU Darmstadt \\
Fachbereich Mathematik \\
Schlo\ss gartenstra\ss e 7 \\
64289 Darmstadt \\
Germany \\
e-mail: {\tt gramlich@mathematik.tu-darmstadt.de} \\
{\tt mars@mathematik.tu-darmstadt.de}

\vspace{.8cm}

\noindent First author's alternative address:

\noindent University of Birmingham \\
School of Mathematics \\
Edgbaston \\
Birmingham \\
B15 2TT \\
United Kingdom \\
e-mail: {\tt ralfg@maths.bham.ac.uk}

\end{document}